\documentclass[11pt]{amsart}
\usepackage{amssymb, graphics}
\usepackage{hyperref}
\newtheorem{lemma}{Lemma}
\newtheorem{prop}{Proposition}
\newtheorem{thm}{Theorem}
\newtheorem{cor}{Corollary}
\theoremstyle{definition}
\newtheorem{rem}{Remark}
\newtheorem{defn}{Definition}
\newtheorem{ex}{Example}

\newcounter{numl}

\newcommand{\labelnuml}{\textup{(\roman{numl})}}
\newenvironment{numlist}{\begin{list}{\labelnuml}%
{\usecounter{numl}\setlength{\leftmargin}{0pt}%
\setlength{\itemindent}{2\parindent}%
\setlength{\itemsep}{\smallskipamount}\def
\makelabel ##1{\hss \llap {\upshape ##1}}}}{\end{list}}

\newcommand{\R}{{\mathbb R}}
\newcommand{\C}{{\mathbb C}}

\newcommand{\cO}{{\mathcal O}}

\newcommand{\symprod}{\mathbin{\raise1pt\hbox{$\scriptstyle\bigcirc$}}}

\newcommand{\JJ}{\mathcal J}
\newcommand{\RR}{\mathbb R}
\newcommand{\delbar}{\overline{\partial}}
\newcommand{\del}{\partial}
\newcommand{\im}{\mathrm{im}}

\begin{document}

\title[Generalized K\"ahler manifolds with split tangent bundle]
{Generalized K\"ahler manifolds with split tangent bundle}

\author[V. Apostolov]{Vestislav Apostolov}
\address{Vestislav Apostolov \\ D{\'e}partement de Math{\'e}matiques\\
UQAM\\ C.P. 8888 \\ Succ. Centre-ville \\ Montr{\'e}al (Qu{\'e}bec) \\
H3C 3P8 \\ Canada}
\email{apostolo@math.uqam.ca}
\author[M. Gualtieri]{Marco Gualtieri}
\address{Marco Gualtieri\\ M.I.T. Department of Mathematics \\  77 Massachusetts Avenue \\Cambridge MA 02139-4307\\ USA }
\email{mgualt@math.mit.edu}

\thanks{We would like to thank P. Gauduchon, G. Grantcharov and N.~J.~Hitchin for their help and stimulating discussions.}
\date{May 22, 2006}

\begin{abstract}
We study generalized K\"ahler manifolds for which the corresponding
complex structures commute and classify completely the compact
generalized K\"ahler four-manifolds for which the induced complex
structures yield opposite orientations.
\end{abstract}

\maketitle

\section{Introduction}
The notion of a {\it generalized K\"ahler structure} was introduced
and studied by the second author in \cite{gualtieri}, in the context
of the theory of generalized geometric structures initiated by
Hitchin in \cite{hitchin}. Recall that a generalized K\"ahler
structure is a pair of commuting complex structures $(\JJ_1,\JJ_2)$
on the vector bundle $TM \oplus T^*M$ over the smooth manifold
$M^{2m}$, which are:
\begin{enumerate}
\item[$\bullet$] integrable with respect to the (twisted) Courant bracket on $TM\oplus T^*M$,
\item[$\bullet$] compatible with the natural inner-product $\langle \cdot, \cdot \rangle$ of signature $(2m,2m)$ on $TM \oplus T^*M$,
\item[$\bullet$] and such that the quadratic form $\langle {\mathcal J}_1\ \cdot,\JJ_2\ \cdot \rangle$ is definite on $TM \oplus T^*M$.
\end{enumerate}

It turns out~\cite{gualtieri} that such a structure on $TM \oplus
T^*M$ is equivalent to a triple $(g,J_+,J_-)$ consisting of a
Riemannian metric $g$ and two integrable almost complex structures
$J_{\pm}$ compatible with $g$, satisfying the integrability relations
$$d^c_+ F_+ +d^c_- F_- = 0, \ dd^c_\pm F_\pm = 0,$$
where $F_{\pm} = gJ_{\pm}$ are the fundamental 2-forms of the
Hermitian structures $(g,J_{\pm})$, and $d^c_\pm$ are the
$i(\delbar-\del)$ operators associated to the complex structures
$J_\pm$.  The closed 3-form $H=d^c_+ F_+=-d^c_-F_-$ is called the
\emph{torsion} of the generalized K\"ahler structure.

These conditions on a pair of Hermitian structures were first
described by Gates, Hull and Ro\v{c}ek~\cite{physicists} as the
general target space geometry for a $(2,2)$ supersymmetric sigma
model.

As a trivial example we can take a K\"ahler structure $(g,J)$ on $M$
and put $J_+=J$, $J_-= \pm J$ to obtain a solution of the above
equations. One can ask, more generally, the following

\vspace{0.2cm} \noindent {\it Question 1.} {When does a compact
complex manifold $(M,J)$ admit a generalized K\"ahler structure
$(g,J_+,J_-)$ with $J=J_+$?}

\vspace{0.2cm} The case of interest is when $J_+ \neq \pm J_-$,
i.e.~when the generalized K\"ahler structure does not come from a
genuine K\"ahler structure on $(M,J)$. In this paper, we  refer to
such generalized K\"ahler structures as {\it non-trivial}.

Despite a growing number of explicit
constructions~\cite{AGG,gualtieri-et-al,hitchin2, kobak,Lin-Tolman},
the general existence problem for non-trivial generalized K\"ahler
structures remains open.  On the other hand, there are a number of
known obstructions, or conditions that the existence of a
generalized K\"ahler structure imposes on the underlying complex
manifold, which we now describe.

Firstly, it follows from the definition that for a complex manifold
$(M,J)$ to admit a compatible generalized K\"ahler structure it must
also admit a Hermitian metric whose fundamental 2-form is $\partial
{\bar \partial}$-closed. This condition on $(M,J)$ is familiar in
Hermitian geometry. It is trivially satisfied if $(M,J)$ is of
K\"ahler type (i.e.~$(M,J)$ admits a K\"ahler metric). When $M$ is
compact and four-dimensional ($m=2$), a result of Gauduchon~
\cite{gauduchon1} affirms that any Hermitian conformal class
contains a metric with $\partial {\bar \partial}$-closed fundamental
form. Hermitian metrics with $\partial{\bar
\partial}$-closed fundamental form naturally appear in  the study of
local index theory~\cite{bismut}, on the moduli space of stable
vector bundles~\cite{lubke-teleman}, and have been much discussed in
the physics literature where they are referred to as `strong
K\"ahler with torsion' structures. Complex manifolds admitting such
Hermitian metrics are the subject of a number of other interesting
results~\cite{egidi,fino-parton-salamon,granch-granch-poon,kutche,spindel-et-al}.
Examples from \cite{egidi}, together with the results of
\cite{fino-granch} and \cite{fino-parton-salamon}, show that there
are compact complex manifolds of any dimension $2m>4$ which do not
admit any Hermitian metric with $\partial {\bar
\partial}$-closed fundamental form.

Secondly, Hitchin~\cite{hitchin2} showed that if $(M,J)$ carries a
generalized K\"ahler structure $(g,J_+, J_-,H)$ such that $J=J_+$
and $J_{+},J_-$ do not commute, then the commutator defines a {\it
holomorphic Poisson structure} $\pi = [J_+,J_-]g^{-1}$ on $(M,J)$.
In the case when $H^0(M, \wedge^2 TM) =0$, for instance, this
result implies that for any compatible generalized K\"ahler
structure on $(M,J)$, the complex structures $J_{+}$ and $J_-$ must
commute, i.e. $J_+J _- = J_-J_+$.

\vspace{0.2cm} Thus motivated, we study in this paper non-trivial
generalized K\"ahler structures $(g,J_+, J_-)$ for which  $J_+$ and
$J_-$ commute. In this case $Q=J_+J_-$ is an involution of the
tangent bundle $TM$, and thus gives rise to a splitting $TM = T_-M
\oplus T_+M$ as a direct sum of the $(\pm 1)$-eigenspaces of $Q$.
Our first result, Theorem~\ref{main}, proves an assertion first made
in \cite{physicists}, which can be stated as follows: {\it the
sub-bundles $T_{\pm}M$ are tangent to the leaves of two transversal
holomorphic foliations ${\mathcal F}_{\pm}$ on $(M,J_+)$ and $g$
restricts to each leaf to define a K\"ahler metric}.

The fact that $T_{\pm}M$ are both {\it holomorphic} and {\it
integrable} sub-bundles of $TM$ directly relates our existence
problem to a conjecture by Beauville~\cite{beauville}, which states
that the holomorphic  tangent bundle $TM$ of a compact complex
manifold  $(M,J)$ of  K\"ahler type  splits as the direct sum of two
holomorphic integrable sub-bundles if and only if $M$ is covered by
the product of two complex manifolds $M_1\times M_2$ on which the
fundamental group of $M$ acts {\it diagonally}. This conjecture has
been confirmed in various
cases~\cite{beauville,campana-peternell,druel}. Combined with
Hitchin's result~\cite{hitchin2} mentioned above, we obtain a wealth
of  K\"ahler complex manifolds which do not admit non-trivial
twisted generalized K\"ahler structures at all. As pointed out
in~\cite{hitchin3}, such examples include (locally) deRham
irreducible compact K\"ahler--Einstein manifolds with $c_1(M)<0$
(see Theorem~\ref{KE-bis} below).

The existence of non-trivial generalized K\"ahler structures for
which $J_+$ and $J_-$ commute thus reduces to the following
question:

\vspace{0.2cm} \noindent {\it Question 2.} Let $(M,J)$ be a compact
complex manifold whose holomorphic tangent bundle splits as a direct
sum of two holomorphic, integrable sub-bundles $T_{\pm} M$. Define a
second almost complex structure $J_-$ on $M$ to be equal to $J$ on
$T_-M$ and to $-J$ on $T_+M$. Does there exist a Riemannian metric
$g$ on $M$ which is compatible with $J_+:=J$ and $J_-$, and such
that $(g,J_{\pm})$ is a generalized K\"ahler structure on $M$?

\vspace{0.1cm} We note that the almost complex structure  $J_-$
defined as above is automatically integrable and commutes with
$J_+$.

The fact that any maximal integral submanifold of $T_{\pm}M$ must be
K\"ahler with respect to a compatible generalized K\"ahler metric
quickly leads to non-K\"ahler examples where the answer to
Question~2  is negative (see Example~\ref{high-dim-ex}). Another
obstruction comes from the fact that the fundamental 2-form of a
compatible generalized K\"ahler metric must be $\partial {\bar
\partial}$-closed (see Example~\ref{solvmanifold}). We are thus led
to suspect that the above existence problem should be more tractable
when $(M,J)$ is of K\"ahler type, and we conjecture that in this
case the answer to our Question~2 is `{\it yes}'. We are able to
establish this in two special cases treated by Beauville
in~\cite{beauville}, namely when $(M,J)$ admits a K\"ahler--Einstein
metric (Theorem~\ref{KE}), and when $(M,J)$ is four-dimensional
($m=2$).

When $M$ is four dimensional, our resuts are much sharper. In this
case there are two classes of generalized K\"ahler structures,
according to whether  $J_{+}$ and $J_{-}$ induce the same or
different orientations on $M$.  In this paper we shall refer to
these cases as generalized K\"ahler structures of {\it bihermitian}
or {\it ambihermitian} type, respectively, though in the terminology
of \cite{gualtieri} they would correspond to generalized K\"ahler
structures of purely even and purely odd type, respectively. Note
that generalized K\"ahler structures of {ambihermitian} type are
precisely those for which $J_+$ and $J_-$ commute and $J_+ \neq \pm
J_-$.

In section \ref{four}, we solve completely the existence problem of
generalized K\"ahler 4-manifolds of ambihermitian type, by proving
the following result.

\begin{thm}\label{main}
A compact complex surface $(M,J)$ admits a generalized K\"ahler
structure of ambihermitian type $(g,J_+,J_-)$ with $J_+ = J$  if and only if the
holomorphic tangent bundle of $(M,J)$ splits as a direct sum of two
holomorpic sub-bundles. Such a complex surface $(M,J)$ is
biholomorphic to one of the following:
\begin{enumerate}
\item[\rm (a)] a geometrically ruled complex surface which is the projectivization of a projectively flat holomorphic vector bundle over a compact Riemann surface;
\item[\rm (b)] a {bi-elliptic} complex surface, i.e.~a complex surface finitely covered by a complex torus;
\item[\rm (c)] a compact complex surface of Kodaira dimension $1$ and even first Betti number, which is an elliptic fibration over a compact Riemann surface, whose only singular fibres are multiple smooth elliptic curves;
\item[\rm (d)] a compact complex surface of general type, uniformized by the product of two hyperbolic planes ${\mathbb H}\times {\mathbb H}$ and with fundamental group acting diagonally on the factors.
\item[\rm (e)] A Hopf surface, with universal covering space $\C^2 \setminus
\{(0,0)\}$ and fundamental group generated by a diagonal
automorphism $(z_1,z_2) \mapsto (\alpha z_1, \beta z_2)$ with
$0<|\alpha| \le  |\beta|<1$, and a diagonal automorphism $(z_1,z_2)
\mapsto (\lambda z_1, \mu z_2)$ with $\lambda, \mu$ primitive
$\ell$-th roots of $1$.
\item[\rm(f)] An Inoue surface in the family $S_{\mathcal M}$ constructed in \cite{inoue}.
\end{enumerate}
On any of the above complex surfaces there exists a family {\rm
(}depending on one arbitrary smooth function on $M${\rm )} of
generalized K\"ahler structures of ambihermitian type.
\end{thm}
To prove this theorem we use the fact that the commuting complex
structures give rise to a splitting of the holomorphic tangent
bundle of $(M,J_+)$ into two holomorphic line bundles $T_{\pm} M$.
Using this splitting and the methods of \cite{gauduchon1},
we describe the set of all generalized K\"ahler structures of
ambihermitian type on such a complex surface. We thus establish a
one-to-one correspondence between four-manifolds admitting
generalized K\"ahler structures of ambihermitian type and complex
surfaces with  split holomorphic tangent
bundle. The latter class of complex surfaces has been studied by
Beauville~\cite{beauville}. We use his classification and some
results from \cite{wall} to derive Theorem~\ref{main}.

We further refine our classification by considering the
\emph{untwisted} case, i.e. when $[H]=0\in H^3(M,\RR)$,  and the
\emph{twisted} case, where $[H]$ is nonzero. We show, by using the
fundamental results of Gauduchon~\cite{gauduchon1, gauduchon2}, that
untwisted generalized K\"ahler structures on compact four-manifolds
can only exist when the first Betti number is even; likewise in the
\emph{twisted} case,  any generalized K\"ahler 4-manifold must have
odd first Betti number (Corollary~\ref{B1}).

\section{Hermitian geometry}
In this section we present certain key properties of Hermitian
manifolds which we will need in the later sections, giving special
attention to the four-dimensional case. Let $M$ be an oriented
$2m$-dimensional manifold. A {\it Hermitian structure} on $M$ is
defined by a pair $(g,J)$ consisting of a Riemannian metric $g$ and
an integrable almost complex structure $J$, which are {\it
compatible} in the sense that $g(J\cdot, J\cdot) = g(\cdot, \cdot)$.
The Hermitian structure $(g,J)$ is called {\it positive} if $J$
induces the given orientation on $M$ and {\it negative} otherwise.

The complex structure $J$ induces a decomposition
$TM\otimes\C=T^{1,0}M\oplus T^{0,1}M$ of the complexified vectors
into $\pm i$ eigenspaces, and hence defines the usual bi-grading of
complex differential forms
$$\Omega^k(M) \otimes \C = \bigoplus_{p+q=k} \Omega^{p,q}(M),$$
where we let $J$ act on $T^*M$ by $(J\alpha)(X)=-\alpha(JX)$, so
that it commutes with the Riemannian duality between vectors and
1-forms: $(J\alpha)^{\sharp} = J \alpha^{\sharp}$.

The product structure $\wedge^2 J$ induces a splitting of the real
2-forms into $\pm 1$ eigenspaces:
$$\Omega^2(M)=\Omega^{J,+}(M)\oplus\Omega^{J,-}(M),$$
whose complexification is simply $\Omega^{J,+}(M)\otimes\C =
\Omega^{1,1}(M)$ and $\Omega^{J,-}(M)\otimes\C =
\Omega^{2,0}(M)\oplus\Omega^{0,2}(M)$.  Furthermore, the {\it fundamental
2-form} $F=gJ$, a real $(1,1)$-form of square-norm $m$, defines a
$g$-orthogonal splitting $\Omega^{J,+}(M)={\mathcal C}^{\infty}(M)\cdot F \oplus
\Omega^{J,+}_0(M)$.  In this way we obtain the $U(m)$ irreducible
decomposition of real 2-forms:
$$\Omega^2(M) = {\mathcal C}^{\infty}(M)\cdot F \oplus \Omega^{J,+}_0(M) \oplus \Omega^{J,-}(M).$$

On a positive Hermitian 4-manifold, the above $U(2)$ splitting of
$\Omega^2(M)$ is compatible with the $SO(4)$ decomposition
$\Omega^2(M)=\Omega^+(M)\oplus\Omega^-(M)$ into self-dual and anti-self-dual
forms, as follows:
\begin{equation}\label{U(2)-split}
\Omega^+(M)= {\mathcal C}^{\infty}(M) \cdot F \oplus \Omega^{J,-}(M); \ \ \Omega^-(M) =
\Omega_0^{J,+}(M).
\end{equation}
For a negative Hermitian structure the r\^oles of $\Omega^+(M)$ and
$\Omega^-(M)$ in the above identifications are interchanged.  Thus, on
an oriented Riemannian four-manifold $(M,g)$, we obtain the
well-known correspondence between smooth sections in $\Omega^+(M)$
(resp. $\Omega^-(M)$) of square-norm 2 and positive (resp. negative)
almost Hermitian structures $(g,J)$. Whereas the existence of such
smooth sections is a purely topological problem, the existence of
integrable ones depends essentially on $g$. This is measured (at
least at a first approximation) by the structure of the Weyl
curvature tensor $W$, cf.~\cite{AG,pontecorvo,salamon}.

The {\it Lee form} $\theta\in\Omega^1(M)$ of a Hermitian structure
is defined by
\begin{equation}\label{Lee0}
dF\wedge F^{m-2} = \frac{1}{(m-1)}\theta \wedge F^{m-1},
\end{equation}
or equivalently $\theta = J \delta^g F$ where $\delta^g$ is the
co-differential with respect to the Levi--Civita connection $D^g$ of
$g$. Since $J$ is integrable, $dF$ measures the deviation of $(g,J)$
from a K\"ahler structure (for which $J$ and $F$ are parallel with
respect to $D^g$). We have the following expression for $D^gF$ (see
e.g.~\cite[p.148]{KN}):
\begin{equation}\label{Kob-Nom}\begin{split}
2g((D^g_X J)Y,Z) 
&= d^cF({X,Y,JZ}) + d^cF(X,JY,Z),
\end{split}
\end{equation}
where $d^c=i(\bar\partial - \partial)$, so that $d^cF= \wedge^3 J(
dF)$ is a real 3-form of type $(1,2)+(2,1)$.

In four dimensions, \eqref{Lee0} reads as
\begin{equation}\label{Lee}
dF = \theta \wedge F,
\end{equation}
and \eqref{Kob-Nom} becomes (see e.g. \cite{gauduchon1, vaisman1})
\begin{equation}\label{DF}
D^g_X F = \tfrac{1}{2}( X^{\flat}\wedge J\theta + JX^{\flat}\wedge
\theta),
\end{equation}
where $X^{\flat}=g(X)$ denotes the $g$-dual 1-form to $X$.  We see
from this that a Hermitian 4-manifold is K\"ahler if and only if
$\theta=0$.

The existence of a K\"ahler metric on a compact complex manifold
$(M^{2m},J)$ implies the Hodge decomposition of the de~Rham
cohomology groups
$$H^k_{dR}(M,\C)\cong \bigoplus_{p+q=k} H^{p,q}_{\bar \partial}(M),$$
where $H^{p,q}_{\bar \partial}(M)$ denote the Dolbeault cohomology
groups.  This, together with the equality $H^{p,q}_{\bar \partial}
(M)\cong {\overline {H^{q,p}_{\bar\partial} (M)}},$  implies that
the odd Betti numbers of a complex manifold admitting a K\"ahler
metric must be even. When $m=2$, it turns out that this condition is
also sufficient.

\begin{thm}\label{kahlerian} \cite{buchdahl,lamari,siu,todorov}
Let $M$ be a compact four-manifold endowed with an integrable almost
complex structure $J$. Then there exists a compatible K\"ahler
metric on $(M,J)$ if and only if $b_1(M)$ is even.
\end{thm}
This important result was first established by
Todorov~\cite{todorov} and Siu~\cite{siu}, using the Kodaira
classification of compact complex surfaces. Direct proofs were found
recently by Buchdahl~\cite{buchdahl} and Lamari~\cite{lamari}.

Since we deal with complex manifolds of non-K\"ahler type (i.e.~do
not admit any K\"ahler metric), we recall the definition of the
$\partial {\bar \partial}$-cohomology groups:
$$H^{p,q}_{\partial {\bar \partial}}(M) := \{d\textrm{-closed} \ (p,q)\textrm{-forms}\}/\partial{\bar \partial}\{(p-1,q-1)\textrm{-forms}\}.$$
Note that there is a natural map
$$\iota : H^{p,q}_{\partial {\bar \partial}}(M) \to H^{p,q}_{\bar \partial}(M).$$
When $(M,J)$ is of K\"ahler type, the well-known $\partial {\bar
\partial}$-lemma (see e.g.~\cite{demailly}) states that the above
map is in fact an isomorphism:
\begin{prop}\label{dd^c-lemma}{\rm ($\partial {\bar \partial}$-lemma)}
If $(M,J)$ is a compact complex manifold admitting a K\"ahler
metric, then $\iota : H^{p,q}_{\partial {\bar \partial}}(M)\to
H^{p,q}_{\bar
\partial}(M)$ is an isomorphism.
\end{prop}
The $\partial {\bar \partial}$-lemma also holds on some non-K\"ahler
manifolds, for example on all non-projective Moi\u{s}ezon manifolds.
In fact, the $\partial {\bar \partial}$-lemma is preserved under
bimeromorphic  transformations and, therefore, holds on any compact
complex manifold which is bimeromorphic to a K\"ahler manifold
(i.e.~is in the so-called {\it Fujiki class ${\mathcal C}$}),
cf.~\cite{demailly}.

\vspace{0.2cm} While the existence of K\"ahler metrics on a compact
complex manifold $(M,J)$ is generally obstructed, a fundamental
result of Gauduchon~\cite{gauduchon1} states that on any {compact}
conformal Hermitian manifold $(M,c,J)$, there exists a unique (up to
scale) Hermitian metric $g \in c$, such that its Lee form $\theta$
is co-closed, i.e. satisfies $\delta^g \theta =0$.  Such a metric is
called a {\it standard} metric of $c$. By \eqref{Lee0}, a standard
metric of $(c,J)$ can be equivalently defined by the equation
$$2i\partial {\bar \partial} F ^{m-1}= dd^c (F^{m-1})=0. $$

We now recall how, in four dimensions, the harmonic properties of
the Lee form with respect to a standard metric are related the
parity of the first  Betti number (compare with
Theorem~\ref{kahlerian} above).
\begin{prop}\label{gauduchon}\cite{gauduchon1, gauduchon2} Let $M$ be a compact four-manifold endowed with a conformal class $c$ of Hermitian metrics, with respect to an integrable almost complex structure $J$. Let $g$ be a standard Hermitian metric in $c$. Then the following two conditions are equivalent:
\begin{enumerate}
\item[\rm (i)] The first Betti number $b_1(M)$ is even.
\item[\rm (ii)] The Lee form $\theta$ of $g$ is co-exact.
\end{enumerate}
\end{prop}
\begin{proof}
For the sake of completeness we outline a proof of this result. Let
$M$ be a compact four-manifold endowed with a standard Hermitian
structure $(g,J)$, and $F$ and $\theta = J \delta^g F$ be the
corresponding  fundamental 2-form and Lee 1-form (with $\delta^g
\theta =0$).

\vspace{0.2cm} We first prove that if $b_1(M)$ is even, then
$\theta$ is co-exact (this is \cite[Th\'eor\`eme II.1]{gauduchon1}).
Applying the Hodge $*$ operator to $\theta$, this is equivalent to
showing that $d^cF$ is exact.  Recall that $2i \partial {\bar
\partial}F= dd^c F=0$ because $g$ is standard. By
Theorem~\ref{kahlerian}, there exists a K\"ahler metric on $(M,J)$
and then, by Proposition~\ref{dd^c-lemma},
$${\bar \partial} F = \partial {\bar \partial} \alpha,$$
for some $(0,1)$-form $\alpha =\xi - i J\xi$. We deduce  $d^cF=
dd^c\xi$, as required.

\vspace{0.2cm} In the other direction, we have to prove that if
$\theta$ is co-exact then $b_1(M)$ is even.  We reproduce an
argument from \cite{gauduchon2}.  With respect to a standard metric $g$,
the forms $\theta$ and $J\theta = - \delta^g F$ are both co-closed,
and therefore the $(0,1)$-form $\theta^{ 0,1} : = \theta - i J\theta
$ is ${\bar \partial}$-coclosed. In terms of Hodge decomposition,
this reads as
$$\theta^{0,1} = \theta^{0,1}_h + {\bar \partial} ^* \Phi,$$
where $\Phi\in\Omega^{0,2}(M)$ and $\theta^{0,1}_h$ is the $({\bar
\partial} {\bar \partial}^* + {\bar \partial}^* {\bar
\partial})$-harmonic part of $\theta^{0,1}$. Note that $\Phi =
\alpha + i \beta$ where $\alpha, \beta \in \Omega^{J,-}(M)$ and $\alpha
(\cdot, \cdot) : = - \beta (J\cdot, \cdot)$.

We first claim that if $\theta^{0,1}_h=0$, then $\phi= F + \beta$ is
a harmonic self-dual 2-form. Indeed, since  $J$ is integrable, it
satisfies $(D^g_{JX} J)(JY) = (D^g_X J)(Y)$ (see \eqref{Kob-Nom}),
and therefore $J (\delta^g \beta) = \delta^g \alpha$, i.e.
$$ \theta - iJ\theta = {\bar \partial} ^* \Phi = \delta^g \Phi = \delta^g \alpha + i \delta^g \beta.$$
It follows that $J\theta = - \delta^g \beta$, and thus
$\delta^g \phi = J\theta + \delta^g \beta =0.$

By a well-known result of Kodaira (see e.g.~\cite{bpv}), a compact
complex surface has even $b_1(M)$ if and only if the dimension
$b_+(M)$ of the space of harmonic self-dual 2-forms on $(M,g)$ is
equal to $2h^{2,0}(M) + 1$, where $h^{2,0}(M) = {\rm dim}_{\mathbb
C} H^{2,0}_{\bar \partial}(M)$; otherwise $b_+(M)=2h^{2,0}(M)$. It
follows that $b_1(M)$ is even if and only if $b_+(M) >  2{\rm
dim}_{\mathbb C} H^{2,0}_{\bar \partial}(M)$.

Therefore, it suffices to show that $\theta^{0,1}_h=0$, provided
that $\theta$ is co-exact (because $\phi$ will be then a harmonic
self-dual 2-form which is not a real part of a holomorphic
$(2,0)$-form). To this end, we consider the natural map $\kappa :
H^1_{dR}(M) \to H^{0,1}_{\bar \partial}(M) \cong H^{1}(M, {\cO})$
from  de Rham to Dolbeault cohomology given by $\xi\mapsto
\xi^{0,1}$ on representatives. One easily checks that $\kappa$ is
well-defined and injective. Moreover, by the Noether formula (see
e.g.~\cite{bpv}), $\kappa$ is an isomorphism of (real) vector spaces
if and only if $b_1(M)$ is even; otherwise, the image of
$H^1_{dR}(M)$ in $H^{0,1}_{\bar
\partial} (M)$ is of real codimension one.

For any element $\xi^{0,1}=\xi - iJ\xi$ in the image of $\kappa$, we calculate its $L_2$-hermitian  product with $\theta^{0,1}_h$:
\begin{equation*}
\begin{split}
 \langle \theta_h^{0,1}, \xi^{0,1} \rangle_{L_2} &=  \langle \theta^{0,1}, \xi^{0,1} \rangle_{L_2} - \langle {\bar \partial}^* \Phi, \xi^{0,1} \rangle_{L_2} \\
 &= \langle \theta^{0,1}, \xi^{0,1} \rangle_{L_2} - \langle  \Phi, {\bar \partial} \xi^{0,1} \rangle_{L_2} \\
 &= \langle \theta^{0,1}, \xi^{0,1} \rangle_{L_2} = \frac{1}{2}(\theta, \xi)_{L_2} + \frac{i}{2}(J\theta, \alpha)_{L_2} \\
 &= \frac{1}{2}(\theta, \xi)_{L_2}  -\frac{i}{2}(\delta^gF, \alpha)_{L_2} = \frac{1}{2}(\theta, \xi)_{L_2}.
 \end{split}
 \end{equation*}
It follows that $\langle \theta_h^{0,1}, \xi^{0,1} \rangle_{L_2}
=0$, if  $\theta$ is co-exact (because $\xi$ is closed).  Thus, in this case,  the image of $\kappa$ is contained in the complex subspace of $H^1_{\bar \partial} (M)$ which is orthogonal to $\theta^{0,1}_h$, and therefore would have real codimension at least 2, unless $\theta_{h}^{0,1}=0$.   \end{proof}

Finally, we review some natural connections which are useful in the
Hermitian context.  An integrable almost complex structure $J$
induces a canonical holomorphic structure on the tangent bundle
$TM$, via the Cauchy--Riemann operator which acts on smooth sections
$X$ and $Y$ of $TM$ by
$${\bar \partial }_X Y := \tfrac{1}{2}([X,Y]+J[JX,Y])=- \tfrac{1}{2} J({\mathcal L}_{Y} J)(X).$$
Identifying $TM$ with the complex vector bundle $T^{1,0}M$, this
operator may be viewed as a partial connection and has the
equivalent expression
\begin{equation}\label{cauchy-riemann-complex}
{\bar \partial}_X Y = [X,Y]^{1,0},
\end{equation}
for any complex vector fields $X$ and $Y$ of type $(0,1)$ and
$(1,0)$, respectively.

In a similar way, any $J$-linear connection $\nabla$ determines a
partial connection ${\bar\partial}^{\nabla}$ on $T^{1,0}$ by
projection, or acting on real vector fields by
\begin{equation}\label{Cauchy-Riemann}
\bar \partial^{\nabla}_X Y = \tfrac{1}{2}(\nabla_X Y + J
\nabla_{JX}Y).
\end{equation}
The operators $\bar \partial$ and $\bar \partial^{\nabla}$ have the
same symbol but do not coincide in general. However, it is
well-known that for any Hermitian structure $(g,J)$, there exists a
unique connection $\nabla$, called the {\it Chern connection} of
$(g,J)$,  which preserves both $J$ and $g$, and such that $\bar
\partial^{\nabla} = \bar \partial$. Note that the Chern connection
$\nabla$ has torsion, unless $(g,J)$ is K\"ahler. It is related to
the Levi--Civita connection $D^g$ by (see e.g. \cite{gauduchon3}):
\begin{equation}\label{chern-connection-higher}
g(\nabla_X Y, Z) = g(D^g_X Y, Z) + \tfrac{1}{2} d^cF({X, JY, JZ}).
\end{equation}

In four dimensions, one uses \eqref{Lee} to rewrite
\eqref{chern-connection-higher} in the following form (cf.
\cite{gauduchon1, vaisman1}):
\begin{equation}\label{chern-connection}
\nabla_X - D^g_X  =  \tfrac{1}{2}\big( X^{\flat}\otimes
\theta^{\sharp} - \theta \otimes X + J\theta(X) J\big),
\end{equation}
where $\theta^{\sharp}= g^{-1}(\theta)$ stands for the vector field $g$-dual to $\theta$.

\section{Generalized K\"ahler structures}\label{three}

As described in the introduction, a generalized K\"ahler structure
on a    manifold $M$ consists of a pair $(\JJ_1,\JJ_2)$ of commuting
generalized complex structures such that $\langle {\mathcal J}_1\
\cdot,\JJ_2\ \cdot \rangle$ determines a definite metric on
$TM\oplus T^*M$.  The generalized complex structures $\JJ_1,\JJ_2$
are integrable with respect to the Courant bracket on sections of
$TM\oplus T^*M$, given by
\[
[X+\xi,Y+\eta]_H = [X,Y] + L_X\eta-L_Y\xi
-\tfrac{1}{2}d(i_X\eta-i_Y\xi) + i_Yi_XH,
\]
which depends upon the choice of a closed 3-form $H$, called the
\emph{torsion} or twisting.  The space of 2-forms $b\in \Omega^2(M)$
acts on $TM\oplus T^*M$ by orthogonal transformations via
\[
e^b(X+\xi) = X+\xi + i_Xb,
\]
and this action affects the Courant bracket in the following way
\[
[e^b(W), e^b(Z)]_H = e^b[W,Z]_{H+db}.
\]
So, if $(\JJ_1,\JJ_2)$ is integrable with respect to the $H$-twisted
Courant bracket, then $(e^{-b}\JJ_1e^{b}, e^{-b}\JJ_2e^{b})$ is
integrable for the $(H+db)$-twisted Courant bracket.

A generalized complex structure $\JJ$, because it is orthogonal and
squares to $-1$, lies in the orthogonal Lie algebra, and therefore
may be decomposed according to the splitting
\[
\mathfrak{so}(TM\oplus T^*M) = \wedge^2 TM \oplus \mathrm{End}(TM)
\oplus \wedge^2 T^*M,
\]
or, in block matrix form,
\[
\JJ = \left(%
\begin{array}{cc}
  A & \pi \\
  \sigma & A \\
\end{array}%
\right),
\]
where $\pi$ is a bivector field, $A$ is an endomorphism of $TM$, and
$\sigma$ is a 2-form. Just as for an ordinary complex structure, the
integrability of $\JJ$ may be expressed as the vanishing of a
Nijenhuis tensor $[\JJ,\JJ]=0$ obtained by extending the Courant
bracket. Restricted to $\wedge^2 TM$, this specializes to the usual
Schouten bracket of bivector fields, requiring that $[\pi,\pi]=0$.
This means that $\pi$ is a Poisson structure.

In \cite{gualtieri}, a complete characterization of the components
of the generalized K\"ahler pair $(\JJ_1,\JJ_2)$ was given in terms
of Hermitian geometry, which we now repeat here.
\begin{thm}[\cite{gualtieri}, Theorem 6.37]\label{mthm}
For any generalized K\"ahler structure $(\JJ_1,\JJ_2)$, there exists
a unique 2-form $b$ and Riemannian metric $g$ such that
\[
e^{-b}\JJ_{1,2}e^b=\frac{1}{2}\left(%
\begin{array}{cc}
  J_+\pm J_- & -(F_+^{-1}\mp F_-^{-1}) \\
  F_+\mp F_- & J_+\pm J_- \\
\end{array}%
\right),
\]
where $J_\pm$ are integrable $g$-compatible complex structures and
$F_\pm = gJ_\pm$ satisfy
\begin{equation}\label{integrabGK}
d^c_+F_++d^c_-F_- = 0, \ dd^c_+ F_+ = 0.
\end{equation}
Conversely, any pair of $g$-compatible complex structures satisfying
condition~(\ref{integrabGK}) define a generalized K\"ahler
structure. Note that the pair $(\JJ_1,\JJ_2)$ is integrable with
respect to the $(H-db)$-twisted Courant bracket where
\[
H = d^c_+ F_+.
\]
\end{thm}
An immediate corollary of this result and the preceding discussion
is that the bivector fields
\begin{equation}\label{poisstruct}
\pi_1 = -F_+^{-1}+F_-^{-1},\ \ \ \ \pi_2= -F_+^{-1}-F_-^{-1}
\end{equation}
are both Poisson structures, a fact first derived
in~\cite{LyakhovichZabzine} directly from~\eqref{integrabGK}.

We also see from the theorem that by taking a bi-Hermitian structure
$(g,J_+,J_-)$ such that $J_+=\pm J_-$, one obtains $d^c_+F_+ =
d^c_-F_-$ and therefore \eqref{integrabGK} reduces to $d^c_+F_+ =
0$, which is nothing but the ordinary K\"ahler condition on
$(g,J_+)$.

As mentioned in the introduction, when $m>2$ the second relation
in~\eqref{integrabGK} imposes a nontrivial constraint on the
underlying complex manifolds $(M,J_{\pm})$: they must admit a
(common) Hermitian metric $g$ for which the fundamental 2-forms are
$dd^c$-closed. Furthermore, if the complex manifold $(M,J_+)$
satisfies the $\partial {\bar
\partial}$-lemma (see Proposition~\ref{dd^c-lemma}), then the
torsion $H=d^c_+ F_+$ of any compatible generalized K\"ahler
structure must be exact.
\begin{prop}\label{kahler-type}
Let $(M,J)$ be a compact complex manifold such that $\iota :
H^{1,2}_{\partial {\bar \partial}}(M) \to H^{1,2}_{\bar
\partial}(M)$ is an isomorphism. Then any generalized
K\"ahler structure on $M$ is \emph{untwisted}, i.e. $[H]=0$.
\end{prop}

We now proceed with an investigation of the class of generalized
K\"ahler structures $(g,J_+,J_-)$ for which the pair of complex
structures \emph{commute} but are unequal, i.e. which satisfy
$[J_+,J_-]=0$ and $J_+\neq \pm J_-$.  In the following theorem, we
show that the splitting
$$TM=T_+M\oplus T_-M,$$
determined by the $\pm 1$-eigenbundles of $Q=J_+J_-$, is not only
integrable, i.e. determines two transverse foliations of $M$, but is
also holomorphic with respect to $J_\pm$, and that the leaves of
each foliation inherit a natural K\"ahler structure.

\begin{thm}\label{main2}
Let $(g,J_+,J_-)$ define a generalized K\"ahler structure with
$[J_+,J_-]=0$. Then the $\pm 1$-eigenspaces of $Q=J_+J_-$ define
$g$-orthogonal $J_\pm$-holomorphic foliations on whose leaves $g$
restricts to a K\"ahler metric.
\end{thm}
\begin{proof}
Let $T_\pm M= \ker (Q\mp {\rm id})=\ker (J_+\pm J_-)$. Since $\ker (J_+\pm
J_-)=\im(J_+\mp J_-)$, we see that $T_\pm M$ coincide with the images
of the Poisson structures
\[
\pi_1 = (J_+-J_-)g^{-1},\ \ \ \ \pi_2 = (J_++J_-)g^{-1}
\]
from~\eqref{poisstruct}.  Therefore $T_\pm M$ are integrable
distributions and determine transverse foliations of $M$. Since $Q$
is an orthogonal operator, we see further that the foliations
defined by its $\pm 1$ eigenvalues must be orthogonal with respect
to the metric $g$.

The complex structures induce decompositions $T_+M\otimes\C = A\oplus
\overline{A}$ and $T_-M\otimes\C = B\oplus \overline{B}$, where
\[
A = T^{1,0}_{J_+} M\cap T^{0,1}_{J_-}M,\ \ \ \ B = T^{1,0}_{J_+}M\cap
T^{1,0}_{J_-}M
\]
are themselves integrable since they are intersections of integrable
distributions.  We now show that $A$ is preserved by the
Cauchy-Riemann operator of $J_+$, proving that $T_+M$ is a
$J_+$-holomorphic sub-bundle. Let $X$ be a $(0,1)$-vector field for
$J_+$ and let $Z\in C^\infty(A)$.  Then
\[
\delbar_XZ = [X,Z]^{1,0}.
\]
Since $T^{1,0}_{J_+} M= A\oplus B$, we may project to these two
components:
\[
\delbar_XZ = [X,Z]_A+[X,Z]_B.
\]
To show that $A$ is $J_+$-holomorphic, we must show the vanishing of
the second term, which upon expanding $X=X_{\overline
A}+X_{\overline B}$, reads
\[
[X,Z]_B=[X_{\overline A},Z]_B + [X_{\overline B},Z]_B.
\]
The first term vanishes since $A\oplus\overline{A}=T_+M\otimes\C$ is
involutive, and the second term vanishes since
$A\oplus\overline{B}=T^{0,1}_{J_-}M$ is involutive. Therefore $A$ is
$J_+$-holomorphic. An identical argument proves that $B$ is
$J_+$-holomorphic, and that both $A,B$ are $J_-$-holomorphic, as
required.

To show that $g$ restricts to a K\"ahler metric on the leaves of
$T_\pm M$, observe that since $J_+=J_-$ along the leaves of $T_-M$, we
have upon restriction $d^c_+F_+ = d^c_-F_-$.  Similarly along the
leaves of $T_+M$ we have $J_+=-J_-$, so that upon restriction,
$d^c_+=-d^c_-$ and $F_+=-F_-$, giving again $d^c_+F_+=d^c_-F_-$. But
since the generalized K\"ahler condition forces
$d^c_+F_+=-d^c_-F_-$, we conclude that both $F_\pm$ are closed upon
restriction to the leaves of either foliation, therefore defining
K\"ahler structures there.
\end{proof}

The holomorphicity of the decomposition $TM=T_+M\oplus T_-M$ proven
above together with the condition $d^c_+ F_+ + d^c_-F_-=0$ also  imply that $Q$ is parallel with respect to the Chern
connections $\nabla^\pm$ of $J_\pm$; in other words, for a
generalized K\"ahler structure with $[J_+,J_-]=0$, the Chern
connections $\nabla^\pm$ have holonomy contained in $U(m_+)\times
U(m_-)$ where $\dim_\R T_\pm M= 2m_\pm$.  We now provide an
alternative proof of this fact, avoiding the use of
Theorem~\ref{mthm}.
\begin{prop}\label{secondprof}
Let $(J_+,J_-)$ be a pair of Hermitian complex structures for the
Riemannian metric $g$, such that $d^c_+F_++d^c_-F_-=0$ and
$[J_+,J_-]=0$. Then $Q=J_+J_-$ is covariant constant with respect to
the Chern connections $\nabla^\pm$.
\end{prop}
\begin{proof}
Since $\nabla^+J_+=0$ by definition, it suffices to show that
$\nabla^+J_-=0$.  From Equation~\eqref{Kob-Nom}, we see that
\[
\nabla^+ - \nabla^- = L,
\]
where $L\in \Omega^1(\mathrm{End}(TM))$ is given by
\[
2g(L_XY,Z)= d^c_+F_+(X,J_+Y,J_+Z) - d^c_-F_-(X,J_-Y,J_-Z).
\]
Consequently, $\nabla^+J_- = \nabla^-J_- + [L,J_-]$.  By definition,
$\nabla^-J_-=0$, and expanding the commutator we obtain
\begin{equation}\label{temp2} \begin{split}
2g([L_X,J_-]Y,Z) &= d^c_+ F_+(X,J_+J_-Y,J_+Z) +  d^c_-
F_-({X,Y,J_-Z})  \\  & +d^c_+ F_+\big(X,J_+Y,J_+J_-Z)  +d^c_-
F_-\big(X,J_-Y,Z).
\end{split}
\end{equation}
If $Y$ is taken in $T_+ M$ and $Z$ in $T_-M$, then the terms cancel
since $d^c_+F_++d^c_-F_-=0$.  If $Y,Z\in T_+ M$, then trivially
$g((\nabla^+_X J_-) Y,Z)=g((\nabla^+_XJ_+)Y,Z)=0$ and similarly for
$Y,Z\in T_-M$. Hence $\nabla^+J_-$ must vanish identically.
Similarly, $\nabla^-J_+=0$, proving the result.
\end{proof}

In fact, this proposition provides an alternative proof not only of
the holomorphicity of $T_\pm M$ but also of their integrability, by
observing that since the torsion of $\nabla^+$ vanishes upon
restriction to $T_\pm M$, we have for $Y,Z\in T_+M$ or $T_-M$,
\[
[Y,Z] = \nabla^+_YZ-\nabla^+_ZY,
\]
and since $\nabla^+Q=0$, $T_\pm M$ are involutive for the Lie
bracket. Applying the same argument to $T_-M$, we obtain an
alternative proof of Theorem~\ref{main2}.

\begin{rem} Along the above lines one can establish the following result:  Let $J_+$ and $J_-$ be a pair of commuting almost complex structures on a $2m$-manifold $M$, such that $J_+$ is {\rm integrable}, and  let $T_{\pm} M$ denote the sub-bundles of $TM$ corresponding to $(\pm 1)$-eigenspaces of $Q= J_+  J_-$.   Then any two of the following three conditions imply the third.
\begin{enumerate}
\item[\rm (a)] $T_{\pm}M$ are integrable sub-bundles of $TM$;
\item[\rm (b)] $T_{\pm}M$ are holomorphic sub-bundles of $TM$ with respect to $J_+$;
\item[\rm (c)] $J_-$ is an integrable almost complex structure.
\end{enumerate}

\end{rem}

\vspace{0.2cm} Let us now return to the existence problem. According
to Theorem~\ref{main2}, we must consider  complex manifolds $(M,J)$
whose tangent bundle splits as a direct sum of two integrable,
holomorphic sub-bundles $T_{\pm}M$; the second complex structure
$J_-$ is obtained from $J_+=J$ by composing with $Q$, the product
structure defining $T_\pm M$.  It is then natural to ask whether there
is a Riemannian metric $g$ on $M$ which is compatible with the
commuting pair $(J_+,J_-)$, satisfying the generalized K\"ahler
condition.  (This is Question~2 of the introduction.)

Locally, the answer is always `{\it yes}'.  Indeed, by using complex
coordinates adapted to the transverse foliations, i.e. a
neighborhood $U= V \times W \subset \C^{m_1}\times \C^{m_2}$ such
that $T_{-}U = TV, \ T_+U=TW$, then for any K\"ahler metrics $g_V$
and $g_W$ on $V$ and $W$, the product metric $g_U := g_V \times
g_{W}$ is K\"ahler with respect to both $J_{\pm}$, and
$(g_U,J_{\pm})$ is a generalized K\"ahler structure.

We now show that if there exists one generalized K\"ahler metric $g$
on $(M,J_+,J_-)$, then there is in fact a whole family parametrized
by smooth functions (This is similar to the variation of a K\"ahler
metric by adding $dd^c f$).  This construction is closely related to
the potential theory developed in \cite{physicists, rocek-et-al}. We
will use the integrable decomposition
\[
TM = T_+M\oplus T_-M,
\]
and the associated decomposition $d = \delta_++\delta_-$ of the
exterior derivative (induced by  the `type' decomposition $\wedge^* T^*M = (\wedge^* T_+M^*) \otimes (\wedge^* T_-M^*)$),   so that, defining $\delta^c_\pm
=[J_+,\delta_\pm]$, we have
\begin{equation}\label{dc}
d^c_\pm = \pm\delta^c_+ + \delta^c_-.
\end{equation}

\begin{prop}\label{potential}
Let $(M,g,J_+,J_-)$ be a generalized K\"ahler structure.  Then, for
any smooth function $f\in C^\infty(M,\R)$ and sufficiently small
real parameter $t$, the 2-form
\begin{equation}\label{variat}
{\tilde F}_+ = F_+ + t (\delta_+\delta^c_+ f - \delta_-\delta^c_-f)
\end{equation}
defines a new Riemannian metric $\tilde g = -\tilde F_+ J_+$  which
is compatible with both $J_{\pm}$, and such that $(\tilde g,
J_{\pm})$ defines a generalized K\"ahler structure with unmodified
torsion class $[H]\in H^3(M,\R)$.
\end{prop}
\begin{proof}
The $J_\pm$-invariant 2-form in~\eqref{variat} defines the
$J_-$-fundamental form $\tilde F_- = \tilde g J_-$, or
\[
{\tilde F}_- = F_- + t (-\delta_+\delta^c_+ f -
\delta_-\delta^c_-f).
\]
We now show that $d^c_+\tilde F_+ + d^c_-\tilde F_-=0$, since
\begin{align*}
d^c_+(\delta_+\delta^c_+-\delta_-\delta^c_-) +
d^c_-(-\delta_+\delta^c_+ - \delta_-\delta^c_-) &=
-\delta^c_+\delta_-\delta^c_- +
\delta_-^c\delta_+\delta^c_+\\
&\ \ \ \ +\delta_+^c\delta_-\delta^c_- -
\delta^c_-\delta_+\delta^c_+\\ &= 0.
\end{align*}
Finally, by the identity
\begin{align*}
d^c_+(\delta_+\delta^c_+ - \delta_-\delta^c_-)&=
\delta_-^c\delta_+\delta_+^c-\delta_+^c\delta_-\delta_-^c\\
&=(\delta_++\delta_-)\delta_+^c\delta_-^c\\
&=d \delta_+^c\delta_-^c,
\end{align*}
we see that $d^c_+(\tilde F_+ - F_+)$ is exact, showing that
$[d^c_+F_+]=[d^c_+\tilde F_+]$, completing the proof.
\end{proof}

The following example shows that the global existence question is
more subtle.
\begin{ex}\label{high-dim-ex}
Take the product  $M = M_1 \times M_2$ of two complex manifolds
$(M_1,J_1)$ and $(M_2, J_2)$, where the latter admits no K\"ahler
metrics at all (see Theorem~\ref{kahlerian}) and put $J_{\pm} := J_1
\pm J_2$ on $TM = TM_1 \oplus TM_2$. Then $J_+$ and $J_-$ commute
and induce the obvious holomorphic splitting of $TM$, but they
cannot admit a compatible generalized K\"ahler metric $g$ (see Theorem~\ref{main2}). In fact,
$(M,J_+,J_-)$ cannot admit any  compatible Riemannian metric $g$
with $d^c_+F_+ + d^c_-F_- =0$ (see Proposition~\ref{secondprof}).  Note that while $(M,J)$ admits no
K\"ahler metric,  $M_1$ can be chosen so that $(M,J)$ does admit
Hermitian metrics with $\partial {\bar
\partial}$-closed fundamental forms.
\end{ex}
By contrast, if $(M,J)$ is a complex manifold of K\"ahler type, we
can always find a Riemannian metric compatible with both $J_+$ and
$J_-$ and such that $d^c_+F + d^c_-F_-=0$, as we now show.
\begin{lemma}~\label{kahler-case}
Let $(M,J_+)$ be a complex manifold of K\"ahler type whose tangent
bundle splits as a direct sum of two holomorphic, integrable
sub-bundles $T_{\pm} M$, and let $J_-=-J|_{T_+M}+J|_{T_-M}$. Then
$M$ admits a Riemannian metric $g$, compatible with both $J_+$ and
$J_-$, satisfying $d^c_+F_+ + d^c_-F_-=0$.
\end{lemma}
\begin{proof}
Let $g_0$ be any K\"ahler metric for $(M,J_+)$; since $J_{\pm}$
commute, the $J_-$-averaged Riemannian metric
$$g(\cdot, \cdot) : = \tfrac{1}{2}(g_0(\cdot, \cdot) + g_0(J_- \cdot, J_- \cdot))$$
is compatible with both $J_{\pm}$. We claim that $g$ has the desired
properties.

To see this, decompose the original K\"ahler form $F_0$ according to
the splitting $\wedge^2 T^*M=\wedge^2 T_+^*M\oplus (T_+^*M\otimes
T_-^*M)\oplus \wedge^2 T_-^*M$, yielding
\[
F_0 = F_{++} + F_{+-}+F_{--}.
\]
Then the fundamental forms for $(g,J_\pm)$ are
\[
F_\pm = \pm F_{++} + F_{--},
\]
and using Equation~\eqref{dc} and the fact $dF_0=0$, we obtain
\begin{align*}
d^c_-F_- &= \delta^c_-(-F_{++}) - \delta_+^cF_{--}\\
&=-d^c_+F_+,
\end{align*}
as required.
\end{proof}

Note that in the above Lemma, the commuting bi-Hermitian structure
$(g,J_\pm)$ is not necessarily generalized K\"ahler, because
although $d^c_+F_++d^c_-F_-=0$, it is not necessarily the case that
$dd^c_+F_+=0$.  We now provide an example where this final condition
cannot be fulfilled.

\vspace{0.2cm}
\begin{ex}~\label{solvmanifold}
We elaborate on an example from \cite{debartolomeis-tomassini} of a
compact $6$-dimensional solvmanifold $M$ which does not admit a
K\"ahler structure.  $M$ is obtained as a compact quotient of a
complex 3-dimensional Lie group (biholomorphic to $\C^3$) whose
complex  Lie algebra $\mathfrak{g}$  is generated by the complex
$(1,0)$-forms $\sigma_1, \sigma_2, \sigma_3$, such that
$$d\sigma_1=0, \ d\sigma_2 = \sigma_1 \wedge \sigma_2,  \ d\sigma_3 = - \sigma_1 \wedge \sigma_3.$$
Thus, $\mathfrak{g}$ (and hence $M$) inherits a natural
left-invariant complex structure $J$ with respect to which the
$\sigma_i$ are holomorphic $1$-forms. Note that $(M,J)$ does not
satisfy the $\partial {\bar \partial}$-lemma because $\sigma_2$  and
$\sigma_3$ are holomorphic but not closed.

It is straightforward to check that there are no left-invariant
Hermitian metrics $g$ on $(\mathfrak{g},J)$ such that the condition
$dd^c F =0$ is satisfied. Since the volume form $v = \sigma_1 \wedge
{\overline \sigma_1} \wedge \sigma_2 \wedge {\overline \sigma_2}
\wedge \sigma_3 \wedge {\overline \sigma_3}$ is bi-invariant, a
standard argument~\cite{belgun, fino-granch} shows that $(M,J)$ does
not admit {\it any} Hermitian metrics with $dd^c$-closed fundamental
form. In particular, $(M,J)$ admits no compatible generalized
K\"ahler structures.

However, we can define a second left-invariant complex structure
$J_-$ on $\mathfrak{g}$ (and hence also on $M$) such that
$T^{1,0}_{J_-}M = {\rm span}_{\C}  \{{\overline \sigma_1}, \sigma_2,
\sigma_3\}$, so that $J_+ :=J$ and $J_-$ are both integrable,
commute and define {\it holomorphic} (and therefore integrable)
sub-bundles $T_{\pm}M$. Furthermore, the left-invariant metric $g_0
= \sum_{i=1}^{3} \sigma_i \otimes {\overline \sigma_i}$ on
$\mathfrak{g}$ defines on $M$ a Hermitian metric which is compatible
with both $J_+$ and $J_-$,  and such that  $d^c_+F_+ + d^c_-F_-=0.$
\end{ex}

For a compact complex manifold of K\"ahler type, $(M,J)$, Beauville
conjectures~\cite{beauville} that $TM$ splits as a direct sum of two
holomorphic integrable sub-bundles if and only if $M$ is covered by
the product of two complex manifolds $M_+\times M_-$ on which the
fundamental group of $M$ acts {\it diagonally}, i.e. $\pi_1(M)$ acts
on each $M_\pm$ and its action on the product is the diagonal
action.  In the case when there is a K\"ahler metric on $(M,J)$
whose Levi-Civita connection preserves $T_+M$ and $T_-M$, the
conjecture follows by the de~Rham decomposition theorem. It has also
been confirmed in other
cases~\cite{beauville,campana-peternell,druel}. We mention here the
following partial result.

\begin{thm}\label{KE}\cite{beauville,kobayashi}
Let $(M,J)$ be a compact complex manifold which admits a
K\"ahler--Einstein metric $g$, and whose tangent bundle splits as a
direct sum of two holomorphic sub-bundles $T_{\pm}M$. Then
$T_{\pm}M$ are parallel with respect to the Levi-Civita connection.
In particular,  $g$ is K\"ahler with respect to both $J_+=J$ and
$J_-=-J|_{T_+M}+J|_{T_-M}$, and therefore $(M, J)$ admits
generalized K\"ahler metrics compatible with $J_+$ and $J_-$.
\end{thm}
\begin{proof}
This is a standard Bochner argument. Let $g$ be a K\"ahler--Einstein
metric on $(M,J)$. The vector bundle $E=\mathrm{End}(TM)$ is a
Hermitian holomorphic bundle with unitary connection $D$ induced by
the Levi--Civita connection.  The Ricci endomorphism of $E$ is
defined by
$$K (Q)=  [R, Q],$$
where  $R \in C^\infty(E)$ is the usual Ricci endomorphism of the
tangent bundle and $Q\in C^\infty(E)$.  Since $g$ is
K\"ahler--Einstein, $K\equiv 0$.

A section $Q\in C^\infty(E)$ is holomorphic if and only if $D''Q=0$,
where $D=D'+D''$ is the usual decomposition of $D$ into partial
connections.  The classical {\it Bochner--Kodaira} identity (see
e.g.~\cite{kobayashi-0,demailly}) implies that for any {\it
holomorphic} section $Q$ of $E$,
\begin{equation}\label{bochner}
\int_M ||D'Q||_g^2 v_g = \int_M g(K(Q),Q) v_g= 0.
\end{equation}
Thus, any holomorphic section of $E$ must be parallel.  Applying
this to $Q=J_+J_-$, we see that $T_\pm M$ are parallel for the
Levi--Civita connection. By the de Rham decomposition theorem,
$(M,g,J)$ must be then a local K\"ahler product of two
K\"ahler--Einstein manifolds tangent to $T_{\pm}M$, respectively.
The claim follows.
\end{proof}

To conclude this section, we wish to indicate that the methods of
Theorem~\ref{main2} and Proposition~\ref{secondprof} can be used to
prove {\it non-existence} results as follows.  When $J_+$ and $J_-$
do not commute, a direct computation using \eqref{temp2} shows that
the commutator $P=[J_+,J_-]$ satisfies
$$ \nabla^+_{X} P  +  J_+ (\nabla^+_{J_+X} P) = 0,$$
provided that $d^c_+F_++d^c_-F_-=0$. It follows that for any
generalized K\"ahler structure $(g,J_{\pm})$,  $P$ defines a
$J_\pm$-{\it holomorphic} bivector field $\pi=Pg^{-1}$. This fact
was first established in \cite{AGG} for the case $m=2$, and by
Hitchin~\cite{hitchin2} in general; the latter work also shows that
$P$ defines a $J_\pm$-holomorphic {\it Poisson structure}, a fact
which follows from the fact that $\pi_1,\pi_2$ are Poisson
structures (see Equation~\eqref{poisstruct}). Therefore, if $(M,J)$
does not carry a non-trivial holomorphic Poisson structure (e.g. if
$H^0(M, \wedge^2(TM)) = 0$), then for any generalized K\"ahler
structure $(g,J_{\pm})$ with $J_+=J$, $J_{+}$ and $J_-$ must
commute.  Then, by Theorem~\ref{main2}, non-trivial generalized
K\"ahler structures do not exist unless the holomorphic tangent
bundle of $(M,J)$ splits.  Using results of
\cite{beauville,campana-peternell,druel} one finds a wealth of
projective complex manifolds such that $H^0(M, \wedge^2(TM))=0$ and
$TM$ does not split.  This argument has been used in \cite{hitchin3}
to prove that a locally de Rham irreducible K\"ahler--Einstein
manifold with $c_1(M)<0$  does not admit any non-trivial generalized
K\"ahler structure, thus establishing a partial converse of
Theorem~\ref{KE}.
\begin{thm}\label{KE-bis}\cite{hitchin3}
Let $(M,J)$ be a compact complex manifold of negative first Chern
class. Then it admits a non-trivial generalized K\"ahler structure
$(g,J_+,J_-)$ with $J_+=J$  if and only if the holomorphic tangent
bundle of $(M,J)$ splits. In this case, $J_+$ and $J_-$ commute.
\end{thm}
\begin{proof}
By the Aubin--Yau theorem~\cite{aubin,yau}, $(M,J)$ admits a
K\"ahler--Einstein metric of negative scalar curvature. A standard
Bochner argument shows $H^0(M, \wedge^2(TM)) = 0$. By the preceding
remarks, for any generalized K\"ahler structure $(g,J_+, J_-)$ with
$J_+=J$,  the complex structures must commute and the result follows
from Theorem~ \ref{KE}.
\end{proof}

\section{Generalized K\"ahler four-manifolds}\label{four}
In dimensions divisible by four, generalized K\"ahler structures
fall into two broad classes, defined by whether the complex
structures $\pm J_+$ and $\pm J_-$ induce the same or different
orientations on the manifold.
\begin{defn}
Let $M$ be a manifold of dimension $4k$. A triple $(g,J_+, J_-)$,
consisting of a Riemannian metric $g$ and two $g$-compatible complex
structures $J_\pm$ with $J_+\neq \pm J_-$, is called a
\emph{bihermitian} structure if $J_+$ and $J_-$ induce the same
orientation on $M$; otherwise, it is called \emph{ambihermitian}.
Similarly, an (am)bihermitian conformal structure is a triple
$(c,J_+,J_-)$, where $c=[g]$ is a conformal class of (am)bihermitian
metrics.
\end{defn}
In this section we will concentrate on the 4-dimensional case, where
we have the following characterization of the generalized K\"ahler
condition in terms of the Lee forms $\theta_\pm$.

\begin{prop}\label{b1}
Let $(g,J_{\pm})$ be an {\rm (}am{\rm )}bihermitian structure on a
four-manifold $M$. Then the condition $d^c_+F_++d^c_-F_-=0$ is
equivalent to $\theta_+ + \theta_-=0$  in the bihermitian case, and
to $-\theta_+ +\theta_- =0$ in the ambihermitian case.  The
condition $dd^c_+F_+=0$ means that $g$ is a standard metric, i.e.
$\delta^g\theta_+=0$. The twisting $[H]$ vanishes if and only if
$\theta_+ = \delta^g\alpha$ for $\alpha\in\Omega^2(M)$, i.e. the Lee
form is co-exact.
\end{prop}
\begin{proof}
By \eqref{Lee}, we have $d^c_{\pm}F_{\pm} =
(J_{\pm}\theta_{\pm})\wedge F_{\pm},$ so that
\begin{equation}\label{temp1}
d^c_+F_++d^c_-F_-=(J_+\theta_+) \wedge F_+  +(J_-\theta_-)\wedge
F_-.
\end{equation}
Note that in the bihermitian case $F_+ \wedge F_+=F_-\wedge F_-$ is
twice the volume form $v_g$, whereas in the ambihermitian case $F_+
\wedge F_+ = - F_-\wedge F_- = 2v_g$. Therefore, applying the Hodge
star operator $*$ to \eqref{temp1} and using the fact that
$\delta^g= -* d *$ when acting on 2-forms, we obtain the result.
\end{proof}

As an immediate corollary of this result, together with
Proposition~\ref{gauduchon}, we obtain\footnote{Alternatively, this
result follows from the generalized Hodge decomposition for
generalized K\"ahler structures proven in~\cite{gualtieri-pq}.}:
\begin{cor}\label{B1}
Let $M$ be a generalized K\"ahler 4-manifold. If the torsion class
$[H]\in H^3(M,\R)$ vanishes, then the first Betti number  must be
even (and hence $M$ is of K\"ahler type); if $[H]\neq 0$ then the
first Betti number must be odd.
\end{cor}

Bihermitian complex surfaces were studied
in~\cite{A,AGG,dloussky,kobak,pontecorvo} and classified for even
first Betti number in~\cite{AGG}, where the classification of
Poisson surfaces~\cite{bartocci-macri} is used, and existence is
only partially proven.  In fact,~\cite{AGG} provides enough to show
that in this case, any bihermitian structure is conformal to a
unique generalized K\"ahler structure, up to scale.
\begin{prop}\label{BH}
Let $(c,J_+,J_-)$ be a bihermitian conformal structure on a compact
four-manifold $M$ with $b_1(M)$ even. Then there is a unique (up to
scale) metric $g\in c$ such that $(g,J_+,J_-)$ is generalized
K\"ahler.
\end{prop}
\begin{proof}
By \cite[Lemma~4]{AGG}, any standard metric $g$ of $(c,J_+)$ (which
is unique up to scale~\cite{gauduchon1}) is standard for $(c,J_-)$
as well, and furthermore $\theta_+ + \theta_-=0$. By
Proposition~\ref{b1}, this is equivalent to the generalized K\"ahler
condition.\end{proof}
Some constructions of these bihermitian
structures can be found in \cite{AGG,gualtieri-et-al,
hitchin2,kobak,Lin-Tolman}, and these prove existence on many (but
not all) of these surfaces.

In the case where the first Betti number is odd, bihermitian
structures have been studied in \cite{A, AGG, dloussky, pontecorvo}.
It follows from the results there that $M$ must be a finite quotient
of $(S^1\times S^3) \sharp k {\overline {\C P}}^2, \ k \ge 0$.  It
is no longer true in this case that the standard metric provides a
generalized K\"ahler metric in all cases.  To the best of our
knowledge, the only known examples of generalized K\"ahler
structures on 4-manifolds with $b_1(M)$ odd are given by standard
metrics in the anti-self-dual bihermitian conformal classes
described in~\cite{pontecorvo}.

We now turn to the ambihermitian case, where we establish a complete
classification of generalized K\"ahler structures.  We start with
the following observation.
\begin{lemma}\label{doubly-almost-complex}
Let $M$ be a four-manifold endowed with a pair $(J_+,J_-)$ of almost
complex structures inducing different orientations on $M$. Then, $M$
admits a Riemannian metric compatible with both $J_\pm$ if and only
if $J_+$ and $J_-$ commute. In this case, the tangent bundle splits
\begin{equation}\label{Q-split}
TM = T_{+}M \oplus T_-M
\end{equation}
as an orthogonal direct sum of Hermitian complex line bundles
defined as the $\pm 1$-eigenbundles of $Q=J_+J_-$.
\end{lemma}
\begin{proof}
Let $g$ be a Riemannian metric on $M$, compatible with $J_+$ and
$J_-$. Fix the orientation on $M$ induced by $J_+$. As discussed in
\S~2,  the fundamental 2-forms $F_+$ and $F_-$ are sections of
$\Omega^+(M)$ and $\Omega^-(M)$, respectively.  Since $\Omega^-(M)$
is in the $+1$-eigenspace of $\wedge^2 J_+$, $F_-$ is
$J_+$-invariant.  Hence $J_+$ and $J_-$ commute. The converse is
elementary.
\end{proof}
The proof of the above lemma shows that the existence of commuting
almost complex structures on a four-manifold is a purely topological
problem (in fact, it is equivalent to the existence of a field of
oriented two-planes~\cite{matsushita}).  Note that a similar
existence problem for pairs of {\it integrable} almost complex
structures on $M$ inducing different orientations was raised in
\cite{beauville0}, and has been almost completely solved in
\cite{kotschick2}.

Our next step is to identify the compact complex surfaces $(M,J)$
that admit a generalized K\"ahler metric $(g,J_+,J_-)$ of
ambihermitian type with $J_+ =J$.
\begin{lemma}\label{split}
Let $(g,J_+,J_-)$ be an ambihermitian structure on a four-manifold
$M$ and let $Q=J_+J_-$ be the almost product structure it defines.
Then the Lee forms satisfy $\theta_+=\theta_-$ if and only if $T_\pm
M$ are holomorphic sub-bundles for $J_\pm$, i.e. $\nabla^\pm Q=0$.
Then the standard metric in the conformal class defines a
generalized K\"ahler metric.

As a result, any compact complex surface $(M,J)$ whose tangent
bundle splits as a sum of holomorphic line bundles admits a
compatible generalized K\"ahler metric.
\end{lemma}
\begin{proof}
If $\theta_+=\theta_-$, then by Proposition~\ref{b1}, we have
$d^c_+F_++d^c_-F_-=0$, and so $T_\pm M$ are holomorphic by
Proposition~\ref{secondprof}.

In the other direction, we use Equation~\eqref{chern-connection} and
the fact that $J_-$ is skew-symmetric to express
$$\nabla^+_X J_- = D^g_X J_- - \tfrac{1}{2} ( X^{\flat} \wedge
(J_-\theta_+)^{\sharp}  + (J_- X)^{\flat} \wedge \theta_+^{\sharp}),
$$ where $\alpha \wedge X = \alpha \otimes X - X^{\flat}\otimes
\alpha^{\sharp}$ for $\alpha \in T^*M$ and $X\in TM$. Finally, by
\eqref{DF}, we obtain
\begin{equation}\label{calculation}
\nabla^+_X J_- = \tfrac{1}{2} ( X^{\flat} \wedge J_-(\theta_- -
\theta_+)^{\sharp} + J_-X^{\flat}\wedge
(\theta_--\theta_+)^{\sharp}).
\end{equation}
It is clear from Equation~\eqref{calculation} that $\nabla^+J_-$
(and hence $\nabla^\pm Q$) vanishes if and only if
$\theta_+=\theta_-$, proving the result.

To prove the final statement, we note that any holomorphic
one-dimensional sub-bundle $T_\pm M \subset TM$ is automatically
integrable, and therefore the almost complex structure $J_- =
-J|_{T_+M}+J|_{T_-M}$ is integrable. By  definition, $J_+=J$ and
$J_-$ commute, and $J_\pm$ induce different orientations.  Clearly
there are Riemannian metrics compatible with both $J_{\pm}$. Then we
may apply the first part of the lemma.
\end{proof}

Now we are ready to prove our  classification results for
ambihermitian generalized K\"ahler structures.

\vspace{0.2cm} \noindent {\bf Proof of Theorem~\ref{main}.} Let
$(M,g,J_+,J_-)$ be a compact generalized K\"ahler four-manifold of
ambihermitian type. By Proposition~\ref{b1} and Lemma~\ref{split}, the
holomorphic tangent bundle of $(M,J_+)$ must split as a direct sum
of two holomorphic line bundles $(T_{\pm}M,J_+)$. Complex surfaces with split tangent bundles were studied and essentially classified by Beauville~\cite{beauville}. We use his results to retrieve the list (a)--(f).

When $b_1(M)$ is even,  the cases that occur according to
\cite{beauville} correspond to the surfaces listed in (a)--(d) of
Theorem~\ref{main}, modulo the fact that our description of the
surfaces in (a) is slightly different from the one in \cite[\S
5.5]{beauville}, and that the existence of a splitting of $TM$ on
{\it any} surface in (c) is not addressed in \cite[\S
5.2]{beauville}.

To clarify these points, we notice that in the case of a ruled
surface $M=P(E) \to \Sigma$,  \cite[Thm.C]{beauville} implies that
the universal cover  is the product $\C P^1\times \mathbb U$, where
${\mathbb U}$ is the universal covering space of $\Sigma$, and  the
diagonal action of $\pi_1(M) = \pi_1(\Sigma)$ gives rise to a
$PGL(2,\C)$ representation of $\pi_1(\Sigma)$, i.e. the holomorphic
bundle $E$ is projectively-flat as claimed in (a).

Note that for any an elliptic fibration $f: M \to \Sigma$ as in (c), the base curve $\Sigma$  can be given  the structure of an orbifold with a $2\pi/m_i$ cone point at each point corresponding to a fibre of multiplicity $m_i$ (see,  \cite[\S~5.2]{beauville} and \cite[\S~7] {wall}). Since the Kodaira dimension of $M$ is equal to $1$, the orbifold Euler characteristic of $\Sigma$ must be negative, and therefore $\Sigma$ is a good orbifold uniformized by the hyperbolic space ${\mathbb H}$. Since the first Betti number of $M$ is even, it follows from  \cite[Thm.7.4]{wall} the universal covering space of $M$ is $\C \times
{\mathbb H}$, on which the fundamental group $\pi_1(M)$ acts diagonally by isometries of the canonical product K\"ahler metric.

When $b_1(M)$ is odd,
the possible cases are described in
\cite[\S\S~(5.2),(5.6),(5.7),(5.8)]{beauville}. To prove that the only complex surfaces that really occur are those listed in (e) and (f) in Theorem~\ref{main} we have to exclude the possibility that $(M,J_+)$ is an
elliptic fibration of Kodaira dimension $1$, odd $b_1(M)$,  and with only multiple singular fibres with smooth reduction. It is shown in \cite[\S~(5.2)]{beauville} that for the holomorphic
tangent bundle of such a surface to split,   it must be covered by a
product of simply connected Riemann surfaces on which the
fundamental group acts diagonally. On the other hand, any  elliptic
surface $M$ with Kodaira dimension $1$ and $b_1(M)$ odd is
finitely covered by an elliptic fiber bundle $M'$ over a compact
Riemann surface of genus $>1$,  which has trivial monodromy~(cf.~\cite[p.139]{wall}).  Since $M$ (and hence $M'$) is not K\"ahler, $b_1(M')$ is odd too. Wall \cite[p.141]{wall} showed that the
universal cover of such an $M'$ is $\C \times {\mathbb H}$ on which $\pi_1(M')$ does not act diagonally. It then follows from Beauville's result  cited above that the holomorphic
tangent bundle of $M'$ (and hence of $M$) does not split.

\vspace{0.2cm}

It remains to establish the existence of generalized K\"ahler
metrics on the complex surfaces listed in Theorem~\ref{main}. We know by
Lemma~\ref{split} that there are ambihermitian metrics $(g,J_+,J_-)$
on $M$, compatible with the holomorphic splitting of $TM$, which are
parametrized by the choice of Hermitian metrics on each of the
factors $T_{\pm}M$, or equivalently by two smooth functions on $M$.
For any such metric $g$, we have $\theta_+^g= \theta_-^g$, where
$\theta_{\pm}^g$ are the corresponding Lee forms (see
Lemma~\ref{split}). Let $g_0$ be a standard metric of $([g],J_+)$,
i.e.~a metric in the conformal class $[g]$ such that $\delta^{g_0}
(\theta_+^{g_0})=0$.  Since
$\theta_+^{g_0} = \theta_-^{g_0}$, the triple $(g_0,J_+,J_-)$
defines a generalized K\"ahler structure of ambihermitian type.  Finally, since the
standard metric is unique up to scale in any conformal
class~\cite{gauduchon1}, we eventually obtain a family of
generalized K\"ahler metrics on $M$, which depend on one arbitrary
smooth function, completing the proof.

\begin{rem}
Some Hopf surfaces described in case (e) of Theorem~\ref{main} (e.g.
those with $\alpha = \beta \in \R$ and $\lambda = \mu =1$) admit a
Riemannian metric $g$ compatible with  a pair of hyper-complex
structures, ${\mathcal HC}_+$ and  ${\mathcal HC}_-$, inducing
different orientations on $M$, and such that for any choice $J_+ \in
{\mathcal HC}_+$ and $J_- \in {\mathcal HC}_-$, $(g,J_+,J_-)$  is a
twisted generalized K\"ahler structure of ambihermitian type. Such
Hopf surfaces do also admit an abundance of twisted generalized
K\"ahler structures of bihermitian type~\cite{AGG,pontecorvo}.
\end{rem}

\end{document}